





\magnification 1200

\multiply\abovedisplayskip by 8\divide\abovedisplayskip by 8
\multiply\belowdisplayskip by 8\divide\belowdisplayskip by 8

\font\titre=cmr12 at 16 pt
\font\tenBbb=msbm10
\font\sevenBbb=msbm7
\font\fiveBbb=msbm5
\newfam\Bbbfam
\textfont\Bbbfam=\tenBbb \scriptfont\Bbbfam=\sevenBbb
\scriptscriptfont\Bbbfam=\fiveBbb
\def\Bbb{\fam\Bbbfam\tenBbb}

\def\R{{\Bbb R}}

\font\ncbf=cmssbx10 at 9pt
\def\abstract#1{{\ninepoint{\narrower\smallskip\noindent
 {\ncbf Abstract.} #1\smallskip}\vskip.2truein}}

\def\hangbox to #1 #2{\vskip1pt\hangindent #1\noindent \hbox to #1{#2}$\!\!$}
\def\myitem#1{\hangbox to 35pt {#1\hfill}}
\newskip\ttglue
\def\ninepoint{\def\rm{\fam0\ninerm}
   \textfont0=\ninerm \scriptfont0=\sixrm \scriptscriptfont0=\fiverm
   \textfont1=\ninei  \scriptfont1=\sixi  \scriptscriptfont1=\fivei
   \textfont2=\ninesy  \scriptfont2=\sixsy  \scriptscriptfont2=\fivesy
 \textfont3=\tenex  \scriptfont3=\tenex  \scriptscriptfont3=\tenex
 \textfont\itfam=\nineit  \def\it{\fam\itfam\nineit}
 \textfont\slfam=\ninesl  \def\sl{\fam\slfam\ninesl}
 \textfont\ttfam=\ninett  \def\tt{\fam\ttfam\ninett}
 \textfont\bffam=\ninebf  \scriptfont\bffam=\sixbf
 \scriptscriptfont\bffam=\fivebf  \def\bf{\fam\bffam\ninebf}
 \tt  \ttglue=.5em plus.25em minus.15em
 \normalbaselineskip=11pt
 \setbox\strutbox=\hbox{\vrule height8pt depth3pt width0pt}
 \let\sc=\sevenrm  \let\big=\ninebig \normalbaselines\rm}

 \font\ninerm=cmr9 \font\sixrm=cmr6 \font\fiverm=cmr5
 \font\ninei=cmmi9  \font\sixi=cmmi6   \font\fivei=cmmi5
 \font\ninesy=cmsy9  \font\sixsy=cmsy6 \font\fivesy=cmsy5
 \font\nineit=cmti9  \font\ninesl=cmsl9  \font\ninett=cmtt9
 \font\ninebf=cmbx9  \font\sixbf=cmbx6 \font\fivebf=cmbx5
 \def\ninebig#1{{\hbox{$\textfont0=\tenrm\textfont2=\tensy
 \left#1\vbox to7.25pt{}\right.$}}}

\def\frac#1#2{{#1\over #2}}
\def\ega{\mathop{=}\limits}
\def\defi{\mathop{\rm def}}

\def\exp{\mathop{\rm exp}\nolimits}
\def\egdef{\ega_{\defi}}
\def\up#1{\raise 1.5pt\hbox{#1}}  

\centerline{\titre Optimal Young's inequality and its converse:}
\medskip
\centerline{\titre a simple proof}
\bigskip

\centerline{by}
\smallskip

\centerline{\bf Franck Barthe}
\vskip 2truecm

\abstract{
We give a new proof of the sharp form of
Young's inequality for convolutions, first proved by Beckner [Be] and
Brascamp-Lieb [BL]. The latter also proved a sharp reverse
inequality in the case of exponents less than $1$. Our proof
is simpler and gives Young's inequality and its converse altogether.}
\vskip 2truecm

 The classical convolution inequality of Young asserts that for all
functions $f\in L^p(\R)$ and $g\in L^q(\R)$ we have
$$ \|f \ast g\|_r \le \|f\|_p \, \|g\|_q,$$
where $p,q,r$ are $\ge 1$ and $1/p + 1/q =1 + 1/r$.
This inequality is sharp only when $p$ or $q$ is one.
The best constants in Young's inequality were found
by Beckner [Be], using tensorisation arguments and
rearrangements of functions. In [BL], Brascamp and Lieb
derived them from a more general inequality, which we will
refer to as {\sl the Brascamp-Lieb inequality}; this Brascamp-Lieb
inequality was also successfully applied to several problems in convex
geometry by K.~Ball (see [B] for one example). The expression
of the best constant for Young's inequality
is rather complicated but can be easily memorized via a
simple principle: it is obtained when $f$ and $g$ are Gaussian
functions on the real line, $f(x) = \exp(-p' x^2)$ and
$g(x) = \exp(- q' x^2)$, where~$p'$ is the conjugate exponent of~$p$. This
principle has been largely developed
by Lieb in the more recent paper [Li]; among many other results, this
paper contains a new proof of the Brascamp-Lieb inequality (let us
also mention [Ba] where we give yet another proof).
\medskip

 A reverse form of Young's inequality was found by Leindler
[Le]: for $0<p,q,r \le 1$ and $f,g$ non-negative,
$$\|f \ast g\|_r \ge \|f\|_p \|g\|_q.$$
Again these inequalities are sharp only when $p$ or $q$ is one. The
sharp reverse inequalities were obtained by Brascamp and Lieb
in the same paper. It is also shown in [BL] that the
reverse Young inequalities imply another important inequality,
the inequality of Leindler and Prekopa, a close relative of the
Brunn-Minkowski inequality ([Le], [Pr]). As far as we know, the
proof from [BL] is the only proof available for this sharp reverse
Young inequality; in our opinion, it is both rather
mysterious and complicated, and uses many ingredients:
tensorisation, Schwarz symmetrisation,
Brunn-Minkowski and some not so intuitive phenomenon
for the measure in high dimension. To the contrary, our argument is
elementary and gives a unified treatment of both cases, the Young
inequality and the reverse inequality.
\medskip

 It is well known that tensorisation arguments
allow to deduce the multidimensional case from the one-dimensional
(see [Be] for example): if the best constant is $C$
for the real line, it will be $C^N$ in the case of $\R^N$.
We state now the precise results. For every $t > 0$, we define $t'$
by $1/t + 1/t' = 1$ (notice that $t'$ is negative when $t < 1$).
Let us introduce for every $t > 0$
$$ C_t = \sqrt{t^{1/t} \over |t'|^{1/t'}}\up.$$
The general multi-dimensional result is as follows:
\bigskip
\goodbreak

\noindent{\bf Theorem 1.}
{\sl Let $p,q,r>0$ satisfy $1/p + 1/q =1 + 1/r$, and let
$f \in L^p(\R^N)$ and $g \in L^q (\R^N)$ be non-negative functions.
\smallskip

 If $p,q,r \ge 1$ then
$$ \|f \ast g\|_r \le
\Bigl(\frac{C_p C_q} {C_r} \Bigr)^N  \|f\|_p \|g\|_q. \leqno(1)$$

 If $p,q,r \le 1$ then
$$ \|f \ast g\|_r \ge
\Bigl(\frac{C_p C_q} {C_r} \Bigr)^N  \|f\|_p \|g\|_q. \leqno(2)$$}
\medskip

 It is easy to check that when $N = 1$ and $p,q  \neq 1$, there
is equality in $(1)$ or $(2)$ for the functions
$f(x) = \exp(-|p'| \, x^2)$ and $g(x) = \exp(-|q'| \, x^2)$.
As was said above, it is enough to
prove the inequalities when $N = 1$. We will prove this case in a
modified form (Theorem~2) for which we introduce some notation.
The condition $1/p + 1/q = 1 + 1/r$ is equivalent to the relation $1/p' +
1/q' = 1/r'$ for the conjugates, and $r'$, $p'$ and $q'$ have the same
sign if $p,q,r > 1$ or $p,q,r < 1$. We set
$$c = \sqrt{r'/q'} \hbox{\ and\  } s = \sqrt{r'/p'}.$$
Notice that $c^2 + s^2 = 1$. We also introduce the constant
$$ K(p,q,r) =  \frac{p^\frac{1}{2p} q^\frac{1}{2q}} {r^\frac{1}{2r}}$$
that will appear several times in the rest of this paper.
We can now state an equivalent form of Theorem~1. Indeed, a simple
change of variables shows that the following Theorem~2 is equivalent
to Theorem~1 when $N = 1$, provided $p,q$ and $r$ are different from 1.
\medskip
\goodbreak

\noindent{\bf Theorem 2.} {\sl Let $p, q, r > 0$ satisfy
$1/p + 1/q = 1 + 1/r$ and either $p, q, r > 1$ or $p, q, r < 1$.
Let $c = \sqrt{r'/q'}$, $s = \sqrt{r'/p'}$, and let $f, g$
be non-negative functions in $L^1(\R)$.
\smallskip

 If $p,q,r > 1$ then
$$ \Bigl( \int_{\R} \Bigl( \int_{\R} f^{1/p} \Bigl(cx-sy \Bigr)
g^{1/q} \Bigl(sx+cy \Bigr) \, dx \Bigr)^r \, dy \Bigr)^{1/r}
   \le K(p,q,r)
\Bigl( \int_{\R} f \Bigr)^{1/p} \Bigl( \int_{\R} g \Bigr)^{1/q}.
$$

 If $p,q,r < 1$ then
$$ \Bigl( \int_{\R} \Bigl( \int_{\R} f^{1/p} \Bigl(cx-sy \Bigr)
   g^{1/q} \Bigl(sx+cy \Bigr) \, dx \Bigr)^r \, dy \Bigr)^{1/r}
\ge K(p,q,r)
\Bigl( \int_{\R} f \Bigr)^{1/p} \Bigl( \int_{\R} g \Bigr)^{1/q}.$$
In both cases, there is equality when $f(x) = \exp(-p x^2)$
and $g(x) = \exp(-q x^2)$.}
\bigskip
\goodbreak

 By the monotone convergence theorem, it is enough to prove Theorem~2
for functions on $\R$ that are dominated by some centered Gaussian
function. Next, it suffices to prove it for continuous, positive functions;
indeed, assume that $f \le G$, where $G(x) = M \exp(-\varepsilon x^2)$
is a centered Gaussian function, for some $\varepsilon > 0$ and $M>0$;
if $G_1(x) = \frac{1}{\sqrt{2 \pi}} \, \exp(-x^2/2)$ and
$G_n(x) = n G_1(nx)$ then $f_n = \min (f \ast G_n ,G)$ tends to $f$
in $L^p$-norm for every $p\ge 1$. Each function $f_n$
is continuous, positive (it vanishes at some $x\in\R$ only if $f$
is the zero function in $L^1$). Let $(g_n)$ be an approximating
sequence for $g$, built in the same way. If Theorem~2 holds for $f_n$
and $g_n$ for all $n$, then it is true for $f$ and $g$ by the dominated
convergence theorem: we first pass to the limit for the inside integral
of the expression at the left side of the inequality;
then, the domination condition is satisfied for the
function of $y$ defined by the inside integral, and we can
conclude. Computations are especially nice if we assume that $f(x)
\le M \exp(-p\varepsilon x^2)$, $g(x) \le M \exp(-q\varepsilon x^2)$
for some $\varepsilon > 0$ and $M>0$.
\bigskip
\goodbreak

 We state now a lemma that is the real crux of the matter.
\smallskip

\noindent{\bf Lemma 1.} {\sl Assume that $p,q,r > 1$ and that
$1/p + 1/q = 1 + 1/r$. Let $f$, $g$, $F$ and $G$
be continuous positive functions in $L^1(\R)$,
such that $\int f = \int F$ and $\int g = \int G$. We have
$$ \eqalign{
     \Bigl( \int \Bigl( \int f^{1/p} (cx - sy) g^{1/q} (sx + cy) \, dx
& \Bigr)^r \, dy   \Bigr)^{1/r}      \cr
 \qquad \qquad \le  \int \Bigl( \int
& F^{r/p}(cX - sY)  G^{r/q}(sX + cY) \, dY \Bigr)^{{1/r}} \, dX.}
\leqno(3)$$}
\medskip

 Let us first comment about this Lemma. The two numbers
$r/p$ and $r/q$ are larger than one, and will play the role
of $1/P$ and $1/Q$ for some $P, Q < 1$. Letting also $1/R = r$,
we get $1/P + 1/Q = 1 + 1/R$, so that the right-hand side of
inequality $(3)$ is similar to the left-hand side, but for exponents less
than $1$. An easy computation will convince the reader that the two
sides of $(3)$ are equal when $f(x) = F(x) = \exp( -p x^2)$ and $g(x) =
G(x) = \exp(-q x^2)$. These facts imply that Lemma~1
contains both Young's inequality and the
reverse Young's inequality with optimal constant. Of course, Lemma~1 is
also valid for non-negative $L^1$ functions by approximation, as
explained before.
\medskip

\noindent {\bf Proof of Lemma~1.}
The proof is based on a parametrization of functions which was used in
[HM] and was suggested by Brunn's proof of the
Brunn-Minkowski inequality. We assume that $f$,
$F$, $g$ and $G$ are continuous and positive functions in $L^1(\R)$,
such that $\int f = \int F$ and $\int g = \int G$.
We may also assume that the left-hand integral in $(3)$ is finite (using
monotone convergence).
Since $\int f = \int F$ and $\int g = \int G$, there exist two
functions $u$ and $v$ from $\R$ to $\R$ such that for all $t$
$$ \int_{- \infty}^{u(t)} f = \int_{-\infty}^{t} F
\hbox{\ \ and\ \ }
\int_{- \infty}^{v(t)} g = \int_{-\infty}^{t} G.$$
Since  $f$, $g$, $F$ and $G$ are continuous and
never vanish, $u$ and $v$ are increasing bijections
of $\R$ and are continuously differentiable. For all $t$,
$$ u'(t).f(u(t)) = F(t)  \hbox{\ \ and\ \ } v'(t).g(v(t)) = G(t).
\leqno(4) $$
The mapping $T$ defined by $T(x,y) = (u(x),v(y))$ is a bijection of $\R^2$.
Let $R$ be the rotation
$$\left( \matrix{
c & -s \cr
s & \hphantom{-}c} \right)$$
in $\R^2$. We consider the change of variable $(x,y) = \Theta (X,Y)$ in
$\R^2$ given by the mapping $\Theta =^t \!  \! R\, TR$; this
means that
$$ x = c \, u(cX-sY) + s\, v(sX+cY),
\ \ y = -s\, u(cX-sY) + c\, v(sX+cY).$$
It is clear that $\Theta$ is a differentiable bijection of $\R^2$. Its
jacobian $J\Theta$ at a point $(X,Y)$ is equal to
$$ J\Theta(X,Y) = u'(cX-sY) v'(sX+cY).$$
We want an upper estimate for the integral (finite by assumption)
$$ I = \Bigl( \int \Bigl( \int f^{1/p} (cx - sy)
 g^{1/q} (sx + cy) \, dx \Bigr)^r \, dy   \Bigr)^{1/r}.$$
Using the $(L^r,L^{r'})$-duality, there exists a positive
function $h$ such that $\|h\|_{r'} = 1$ and
$$I = \int \!\!\! \int  f^{1/p} (cx-sy) g^{1/q} (sx+cy) h(y) \, dx \, dy.$$
By the change of variable $(x,y) = \Theta (X,Y)$,
we see that $I$ is equal to
$$\eqalign{
 \int \!\!\! \int f^{1/p}(u(cX - sY)) g^{1/q}(v(sX + cY))
h(-s\, u(& cX - sY) + c\, v(sX + cY)) \cr
       u'(&cX - sY) v'(sX + cY) \, dX  dY.}$$
In order to shorten the formulas, let us write
$$\eqalign{
U &= u(cX-sY),\, V = v(sX+cY)\cr
U'&= u'(cX-sY), \,V' = v'(sX+cY).\cr
}$$
From the relations $(4)$ we get
$$
\eqalign{
I = & \int \!\!\! \int f^{1/p}(u(cX-sY)) g^{1/q} (v(sX+cY)) h(-sU+cV)
U'V' \, dX \, dY \cr
  = &  \int \Bigl( \int  F^{1/p} (cX-sY) G^{1/q} (sX+cY) h(-sU+cV)
(U')^{1/p'} (V')^{1/q'} \, dY \Bigr)\, dX.}$$
Using H\"{o}lder's inequality for the integral in $Y$ with
parameters $r$ and $r'$, we obtain
$$\eqalign{
I \le \int \Bigl( \int F^{r/p} (cX - sY) G^{r/q} (sX +
& cY) \, dY \Bigr)^{1/r} \cr
\Bigl(
& \int h^{r'}(-sU + cV) (U')^{r'/p'} (V')^{r'/q'}
\, dY \Bigr)^{1/r'} \, dX.}$$
Let $\displaystyle H(X) =
\int h^{r'}(-sU + cV) (U')^{r'/p'} (V')^{r'/q'} \, dY $, then
$$ H(X) = \int h^{r'} (a(X,Y)) (u'(cX-sY))^{s^2} (v'(sX+cY))^{c^2} \, dY,$$
where
$$a(X,Y) = -s\,u(cX-sY)+c\,v(sX+cY).$$
We have
$$ \frac{\partial a}{\partial Y} (X,Y)=s^2 u'(cX-sY)+c^2 v'(sX+cY).$$
By the arithmetic-geometric inequality
$(U')^{s^2}(V')^{c^2} \le s^2 U' + c^2 V'$, hence
$$ H(X) \le  \int h^{r'} (a(X,Y)) \frac{\partial a}{\partial Y} (X,Y) \, dY =
\int h^{r'} = 1.$$
This proves that
$$ I \le  \int
\Bigl( \int F^{r/p} (cX - sY)  G^{r/q} (sX + cY) \, dY \Bigr)^{1/r}
\, dX$$
and this ends the proof of Lemma~1.
\bigskip
\goodbreak

\noindent {\bf Proof of Theorem~2.}
 If we apply Lemma~1 with
$$ f(x) = F(x) = \sqrt{p/\pi} \, \exp(-px^2),
\ \ g(x) = G(x) = \sqrt{q/\pi} \, \exp(-qx^2),$$
there is equality in $(3)$ and both members are equal to $K(p,q,r)$.
Applying $(3)$ with any $f$ and $g$ such that $\int f = \int g = 1$
and with the preceding $F$ and $G$ gives
Theorem~2 for $p, q, r > 1$. Suitably read from right to left, inequality
$(3)$ gives Theorem~2 when the indices are less than $1$. Indeed,
let $p_1, q_1, r_1<1$ be such that $1/p_1 + 1/q_1 = 1 + 1/r_1$, and
assume that $\int f = \int F = \int g = \int G = 1$. If we define
the triple $(p, q, r)$ by $p = p_1/r_1$, $q = q_1/r_1$ and $r = 1/r_1$,
then $p, q, r > 1$ and  $1/p + 1/q = 1 + 1/r$. So inequality $(3)$
is valid for this triple. A straightforward computation gives that
$$ c_1 \,  \egdef \, \sqrt{\frac{r_1'}{q_1'}}= 
\sqrt{\frac{r'}{p'}} = s  \qquad \hbox{\ and\ } \qquad
 s_1 \, \egdef \, \sqrt{\frac{r_1'}{p_1'}} = \sqrt{\frac{r'}{q'}} = c.$$
Thus, inequality $(3)$ raised to the power $r$ becomes
$$
\eqalign{
  \int \Bigl( \int f^{r_1/p_1} (s_1 x - c_1 y)  g^{r_1/q_1}
&(c_1 x + s_1 y) \,  dx \Bigr)^{1/r_1} \, dy  \cr
  \qquad \le \Bigl( \int \Bigl( \int F^{1/p_1}
&(s_1 X - c_1 Y) G^{1/q_1} (c_1 X + s_1 Y) \, dY
   \Bigr)^{r_1} \, dX   \Bigr)^{1/r_1} .}$$
This is exactly the reverse version of $(3)$ for $p_1, q_1, r_1 < 1$
applied to the functions $\check{F},G,\check{f},g$ where
$\check{f}(x)=f(-x)$. As before, choosing
$f(x) = \sqrt{p_1/\pi} \, \exp(-p_1 x^2),$ and
$ g(x) = \sqrt{q_1/\pi} \, \exp(-q_1 x^2)$ implies Theorem~2
when the parameters are less than 1.
\bigskip
\goodbreak

By the previous argument, there is equality in Theorem~2 when
$f(x) = \exp(-p x^2),$ and $ g(x) = \exp(-q x^2)$. We prove
now that up to scalar multiplication, translation and dilatation,
this is the only equality case.
\medskip

\noindent{\bf Theorem 3.} {\sl
Let $p, q, r >0$ be such that $1/p + 1/q = 1 + 1/r$ and
either $p, q, r > 1$ or $p, q, r < 1$. Let
$c = \sqrt{r'/q'}$, $s = \sqrt{r'/p'}$, and let $f,g$ be
two non-negative functions in $L^1(\R)$. Then
$$ \Bigl( \int \Bigl( \int f^{1/p} \Bigl(cx - sy \Bigr)
  g^{1/q} \Bigl(s x + c y \Bigr) \, dx \Bigr)^r \, dy \Bigr)^{1/r}
   = K(p,q,r)
\Bigl( \int f \Bigr)^ {1/p} \Bigl( \int g \Bigr)^ {1/q} \leqno(5)$$
if and only if there exist $a,b \ge 0$, $\lambda>0$ and $y,z \in \R$ such
that for all $x$
$$ \eqalign{
 f(x) & =  a \, \exp(-\lambda p(x-y)^2) \cr
 g(x) & =  b \, \exp(-\lambda q(x-z)^2).
} \leqno(6)$$}
\medskip

\noindent {\bf Proof.}
Using a simple change of variables, one can check that functions
of the form $(6)$ satisfy equality $(5)$. We show now that only 
these functions do. We give the proof for $p, q, r > 1$; the other
 case is similar.
Let us assume first that $f$ and $g$
are continuous, positive and satisfy equality $(5)$. We may
assume that $\int f = \int g =1$. If we set
$F(x) = \sqrt{p/\pi} \, \exp{(-px^2)}$ and
$G(x) = \sqrt{q/\pi} \, \exp{(-qx^2)}$, we get equality in $(3)$.
We follow the proof of Lemma~1 step by
step. First, we know here that the integral $I$ is finite
by equality $(5)$. There must be equality everywhere
in the proof of inequality $(3)$ for $f$, $g$, $F$ and $G$.
In particular the equality when the arithmetic-geometric inequality
was applied implies that for all $X, Y$ (with
the notation from the proof of Lemma~1)
$$ u'(cX - sY) = v'(sX + cY).$$
So there exists $\mu > 0$ such that $u' = v' = \mu$.
Therefore $u(t) = \mu(t-x_0)$ for some $x_0$. Formula $(4)$ implies that
$$\mu f(\mu(t - x_0)) = \sqrt{\frac{p}{\pi}} \, \exp{(-pt^2)},$$
so $f$ is Gaussian with variance $\mu/p$. By the same method we show
that $g$ is Gaussian with variance $\mu/q$.
\smallskip

 For general $f$ and $g$, we need the following lemma, which was communicated
to me by K. Ball (the reader will recognize in $(7)$ the form of the
Brascamp-Lieb inequality; the next Lemma tells something about maximizers
for this inequality). We denote by $\langle . , .\rangle$ the scalar product in
$\R^n$.
\smallskip

\noindent{\bf Lemma 2.}  {\sl
Let $m \ge n$ be integers and $\alpha_i > 0$, $u_i \in \R^n$,
$i = 1\ldots m$. Assume that there exists $M>0$
such that for all non-negative integrable functions
$f_i$, $i = 1,\ldots ,m$ on $\R$, one has
$$ \int_{\R^n} \prod_{i=1}^m f_i^{\alpha_i} (\langle x,u_i \rangle)\, dx
\le M
\prod_{i=1}^m \Bigl( \int_{\R} f_i \Bigr)^{\alpha_i} \leqno(7)$$
and assume that $M$ is the smallest possible constant for
which this is true. Let us call {\rm maximizer}
a $m$-tuple $(f_1,\ldots,f_m)$ of non zero functions for which
inequality $(7)$ is an equality.

If  $(f_1,\ldots,f_m)$ and $(g_1,\ldots,g_m)$ are maximizers, then so is
$(f_1\ast g_1,\ldots,f_m\ast g_m)$.}
\medskip

\noindent{\bf Proof.} We may assume that $(f_1,\ldots,f_m)$ and
$(g_1,\ldots,g_m)$ are maximizers
and $\int f_i =\int g_i =1$, for $i=1, \ldots ,m$. We
define two functions $F$, $G$ on $\R^n$ by
$$ F(x) = \prod_{i=1}^{m} f_i^{\alpha_i} (\langle x, u_i \rangle)
 \hbox{\ \ and\ \ }
G(x) = \prod_{i=1}^{m} g_i^{\alpha_i} (\langle x, u_i \rangle).$$
We know that $\int F = \int G = M$. So
$$ \eqalign{
M^2 =& \Bigl(\int F\Bigr)  \Bigl(\int G\Bigr) = \int F \ast G \cr
    =& \int \int \prod_{i=1}^{m} f_i^{\alpha_i} (\langle x-y, u_i \rangle)
\prod_{i=1}^{m} g_i^{\alpha_i} (\langle y, u_i \rangle) \, dy \, dx \cr
    =& \int \Bigl( \int \prod_{i=1}^{m} \left[ f_i (\langle x, u_i
\rangle - \langle y, u_i \rangle)
g_i(\langle y, u_i \rangle) \right]^{\alpha_i} dy \Bigr) dx.}$$
Applying inequality $(7)$ to the functions
$k_i(t) = f_i (\langle x, u_i \rangle - t) g_i(t)$, we get
$$M^2\le M \int \prod_{i=1}^{m} \Bigl(\int f_i (\langle x, u_i
\rangle-t)g_i(t) \,
dt \Bigr)^{\alpha_i} \, dx = M \int
\prod_{i=1}^{m} \Bigl(f_i \ast g_i\Bigr)^{\alpha_i}(\langle x, u_i
\rangle)\, dx.$$
It follows that $(f_1\ast g_1,\ldots,f_m\ast g_m)$ is a maximizer.
\medskip

 Now we can finish the proof of Theorem~3.
The functions $f$ and $g$ satisfy $(5)$. As the functions
$\Gamma_p(x)= \exp{(-px^2)}$ and
$\Gamma_q(x)= \exp{(-qx^2)}$ have the same extremal property,
the preceding lemma implies that $f \ast \Gamma_p$
and $g \ast \Gamma_q$ have it too. But these functions are positive and
continuous; by the previous argument
they are of the form $(6)$. Using the Fourier transform, one obtains
that $f$ and $g$  are of the form $(6)$.
\vskip 1truecm

{\ninepoint
\noindent{\bf Acknowledgement.}
I want to express my gratitude to Prof. B. Maurey for drawing my
attention to this problem and for many fruitful discussions.}
\vskip 0.5 truecm
\goodbreak

\noindent{\bf References}
\medskip

\myitem{[B]} K.~Ball, {\sl Volume ratios and a reverse isoperimetric
inequality,} J. London Math. Soc. {\bf 44} (1991), 351--359.
\smallskip

\myitem{[Ba]} F.~Barthe, {\sl In\'egalit\'es de Brascamp-Lieb et convexit\'e,}
 C. R. Acad. Sc. Paris  (1997), to appear.
\smallskip

\myitem{[Be]}
W.~Beckner, {\sl Inequalities in {F}ourier analysis,}
 {Annals of Math.} {\bf 102} (1975), 159--182.
\smallskip

\myitem{[BL]}
H.J.~Brascamp and E.H.~Lieb,
{\sl  Best constants in {Y}oung's inequality, its converse and its
  generalization to more than three functions,}
{Adv. Math.} {\bf 20} (1976), 151--173.
\smallskip

\myitem{[HM]}
R.~Henstock and A.M.~Macbeath,
{\sl  On the measure of sum sets. ({I}) the theorems of {B}runn,
  {M}inkowski and {L}usternik,}
{Proc. London Math. Soc.} {\bf 3} (1953), 182--194.
\smallskip

\myitem{[Le]}
L.~Leindler,
{\sl  On a certain converse of {H}\" older's inequality,}
{Acta Sci. Math. Szeged}, {\bf 33} (1972), 217--223.
\smallskip

\myitem{[Li]}
E.H.~Lieb, {\sl Gaussian kernels have only gaussian maximizers,}
{Invent. Math.}, {\bf 102} (1990), 179--208.
\smallskip

\myitem{[Pr]} A.~Prekopa,
{\sl On logarithmically concave measures and functions,}
Acta Sci. Math. {\bf 34} (1973), 335--343.

\vskip 2truecm

\hbox{\noindent\hskip 5truecm
\vbox{\advance\hsize by -5truecm
\ninepoint
\centerline{Equipe d'Analyse et de Math\'ematiques Appliqu\'ees}

\centerline{Universit\'e de Marne-la-Vall\'ee, 2 rue de la Butte Verte}

\centerline{93166 Noisy-Le-Grand CEDEX. France}

\centerline{e-mail: barthe@clipper.ens.fr}}}

\bye